\documentclass[12pt,reqno]{amsart}
\usepackage{fullpage}
\usepackage{times}
\newtheorem{theorem}{Theorem}[section]

\newtheorem{conjecture}[theorem]{Conjecture}
\newtheorem{corollary}[theorem]{Corollary}
\newtheorem{lemma}[theorem]{Lemma}
{
\theoremstyle{remark}
\newtheorem{definition}[theorem]{Definition}
\newtheorem{example}[theorem]{Example}
\newtheorem*{remark}{Remark}
}

\renewcommand{\mod}[1]{{\ifmmode\text{\rm\ (mod~$#1$)}\else\discretionary{}{}{\hbox{ }}\rm(mod~$#1$)\fi}}
\newsymbol\dnd 232D

\newcommand{\F}{{\mathbb Z}}
\newcommand{\Fpp}{\F_p^2}

\newcommand{\lin}[1]{\langle#1\rangle}
\newcommand{\orig}{\text{\bf 0}}
\newcommand{\noto}{\setminus\{\orig\}}

\newcommand{\GG}{{\mathbf G}}
\newcommand{\HH}{{\mathbf H}}
\newcommand{\KK}{{\mathbf K}}

\begin{document}

\title{Lower bounds for sumsets of multisets in $\F_p^2$} 
\author{Greg Martin}
\address{Department of Mathematics \\ University of British Columbia \\ Room 121, 1984 Mathematics Road \\ Vancouver, BC, Canada V6T 1Z2}
\email{gerg@math.ubc.ca}
\author{Alexis Peilloux}
\address{Universit\'e Pierre-et-Marie Curie \\ Facult\'e de Math\'ematiques \\ 4 place Jussieu \\75005 \\ Paris, France}
\email{alexis.peilloux@etu.upmc.fr}
\author{Erick B. Wong}
\address{Department of Mathematics \\ University of British Columbia \\ Room 121, 1984 Mathematics Road \\ Vancouver, BC, Canada V6T 1Z2}
\email{erick.wong@alumni.ubc.ca}


\begin{abstract}
The classical Cauchy--Davenport theorem implies the lower bound $n+1$ for the number of distinct subsums that can be formed from a sequence of $n$ elements of the cyclic group $\F_p$ (when $p$ is prime and $n<p$). We generalize this theorem to a conjecture for the minimum number of distinct subsums that can be formed from elements of a multiset in $\F_p^m$; the conjecture is expected to be valid for multisets that are not ``wasteful'' by having too many elements in nontrivial subgroups. We prove this conjecture in $\Fpp$ for multisets of size $p+k$, when $k$ is not too large in terms of~$p$.
\end{abstract}

\maketitle

\section{Introduction}

Determining the number of elements in a particular abelian group that can be written as sums of given sets of elements is a topic that goes back at least two centuries. The most famous result of this type, involving the cyclic group $\F_p$ of prime order $p$, was established by Cauchy in 1813~\cite{Cauchy} and rediscovered by Davenport in 1935~\cite{Dav1,Dav2}:

\begin{lemma}[Cauchy--Davenport Theorem]
Let $A$ and $B$ be subsets of $\F_p$, and define $A+B$ to be the set of all elements of the form $a+b$ with $a\in A$ and $b\in B$. Then $\#(A+B) \ge \min\{ p, \#A + \#B - 1\}$.
\end{lemma}

\noindent The lower bound is easily seen to be best possible by taking $A$ and $B$ to be intervals, for example.  It is also easy to see that the lower bound of $\#A + \#B - 1$ does not hold for general abelian groups $\GG$ (take $A$ and $B$ to be the same nontrivial subgroup of $\GG$). There is, however, a well-known generalization obtained by Kneser in 1953~\cite{Kne}, which we state in a slightly simplified form that will be quite useful for our purposes (see \cite[Theorem 4.1]{Nat} for an elementary proof):

\begin{lemma}[Kneser's Theorem]
Let $A$ and $B$ be subsets of a finite abelian group $\GG$, and let $m$ be the largest cardinality of a proper subgroup of $\GG$.  Then $\#(A+B) \ge \min\{ \#\GG, \#A + \#B - m\}$.
\end{lemma}

Given a sequence $A = (a_1,\dots,a_k)$ of (not necessarily distinct) elements of an abelian group $\GG$, a related result involves its {\em sumset} $\Sigma A$, which is the set of all sums of any number of elements chosen from $A$ (not to be confused with $A+A$, which it contains but usually properly):
\[
\Sigma A = \bigg\{ \sum_{j\in J} a_j \colon J\subseteq \{1,\dots,k\} \bigg\}.
\]
(Note that we allow $J$ to be empty, so that the group's identity element is always an element of $\Sigma A$.) When $\GG=\F_p$, one can prove the following result by writing $\Sigma A = \{0, a_1\} + \cdots + \{ 0,a_k \}$ and applying the Cauchy--Davenport theorem inductively:

\begin{lemma}
\label{CD pre lemma}
Let $A=(a_1,\dots,a_k)$ be a sequence of nonzero elements of $\F_p$. Then $\#\Sigma A \ge \min\{p,k + 1\}$.
\end{lemma}

\noindent This result can also be proved directly by induction on $k$, and in fact such a proof will discover why the order $p$ of the cyclic group must be prime (intuitively, the sequence $A$ could lie completely within a nontrivial subgroup). For a formal proof, see \cite[Lemma 2]{DEF}. Again the lower bound is easily seen to be best possible, by taking $a_1=\cdots=a_k$.

It is a bit misleading to phrase such results in terms of sequences, since the actual order of the elements in the sequence is irrelevant (given that we are considering only abelian groups). We prefer to use {\em multisets}, which are simply sets that are allowed to contain their elements with multiplicity. If we let $m_x$ denote the multiplicity with which the element $x$ occurs in the multiset $A$, then the definition of $\Sigma A$ can be written in the form
\[
\Sigma A = \bigg\{ \sum_{x\in\GG} \delta_x x \colon 0\le \delta_x \le m_x \bigg\},
\]
where $\delta_x x$ denotes the group element $x+\cdots+x$ obtained by adding $\delta_x$ summands all equal to~$x$.

When using multisets, we should choose our notation with care: the hypotheses of such results tend to involve the total number of elements of the multiset $A$ counting multiplicity, while the conclusions involve the number of distinct elements of $\Sigma A$. Consequently, throughout this paper, we use the following notational conventions:
\begin{itemize}
 \item $|S|$ denotes the total number of elements of the multiset $S$, counted with multiplicity;
 \item $\#S$ denotes the number of distinct elements of the multiset $S$, or equivalently the number of elements of $S$ considered as a (mere) set.
\end{itemize}
In this notation, Lemma~\ref{CD pre lemma} can be restated as:

\begin{lemma}
\label{CD lemma}
Let $A$ be a multiset contained in $\F_p$ such that $0\notin A$. Then $\#\Sigma A \ge \min\{p,|A| + 1\}$.
\end{lemma}

\noindent The purpose of this paper is to improve, as far as possible, this lower bound for multisets contained in the larger abelian group $\Fpp$.  We cannot make any progress without some restriction upon our multisets: if a multiset is contained within a nontrivial subgroup of $\Fpp$ (of cardinality $p$), then so is its sumset, in which case the lower bound $\min\{p,|A| + 1\}$ from Lemma~\ref {CD lemma} is the best we can do. Therefore we restrict to the following class of multisets. We use the symbol $\orig=(0,0)$ to denote the identity element of~$\Fpp$.

\begin{definition}
\label{valid definition}
A multiset $A$ contained in $\Fpp$ is called {\em valid} if:
\begin{itemize}
\item $\orig\notin A$; and
\item every nontrivial subgroup contains fewer than $p$ points of $A$, counting multiplicity.
\end{itemize}
\end{definition}

\noindent The exact number $p$ in the second condition has been carefully chosen: any nontrivial subgroup of $\Fpp$ is isomorphic to $\F_p$, and so Lemma~\ref{CD lemma} applies to these nontrivial subgroups. In particular, any multiset $A$ containing $p-1$ nonzero elements of a nontrivial subgroup will automatically have that entire subgroup contained in its sumset $\Sigma A$, so allowing $p$ nonzero elements in a nontrivial subgroup would always be ``wasteful''.

We believe that the following lower bound should hold for sumsets of valid multisets:

\begin{conjecture}
\label{2d conjecture}
Let $A$ be a valid multiset contained in $\Fpp$ such that $p \le |A| \le 2p-3$. Then $\#\Sigma A \ge (|A|+2-p)p$. In other words, if $|A| = p+k$ with $0\le k\le p-3$, then $\#\Sigma A \ge (k+2)p$.
\end{conjecture}

It is easy to see that this conjectured lower bound would be best possible: if $A$ is the multiset that contains the point $(1,0)$ with multiplicity $p-1$ and the point $(0,1)$ with multiplicity $k+1$, then the set $\Sigma A$ is precisely $\big\{ (s,t)\colon s\in\F_p,\, 0\le t\le k+1 \big\}$, which has $(k+2)p$ distinct elements. Conjecture~\ref {2d conjecture} is actually part of a larger assertion (see Conjecture~\ref{any d conjecture}) concerning lower bounds for sumsets in $\F_p^m$.

One of our results completely resolves the first two cases of this conjecture:

\begin{theorem}
\label{conjecture true for k tiny}
Let $p$ be a prime.
\begin{enumerate}
\item If $A$ is any valid multiset contained in $\Fpp$ with $|A| = p$, then $\#\Sigma A \ge 2p$.
\item Suppose that $p\ge5$. If $A$ is any valid multiset contained in $\Fpp$ with $|A| = p+1$, then $\#\Sigma A \ge 3p$.
\end{enumerate}
\end{theorem}

It turns out that proving part (b) of the theorem requires a certain amount of computation for a finite number of primes (see the remarks following the proof of the theorem in Section~\ref{thm proofs section}). Extending the conjecture to larger values of $k$ would require, by our methods, more and more computation to take care of small primes $p$ as $k$ grows. However, we are able to establish the conjecture when $p$ is large enough with respect to $k$, or equivalently when $k$ is small enough with respect to $p$:

\begin{theorem}
\label{conjecture true for k small in terms of p}
Let $p$ be a prime, and let $2\le k\le \sqrt{p/(2\log p+1)}-1$ be an integer. If $A$ is any valid multiset contained in $\Fpp$ with $|A| = p+k$, then $\#\Sigma A \ge (k+2)p$.
\end{theorem}

A contrapositive version of Theorem~\ref{conjecture true for k small in terms of p} is also enlightening:

\begin{corollary}
\label{main corollary}
Let $p$ be a prime, and let $2\le k\le \sqrt{p/(2\log p+1)}-1$ be an integer. Let $A$ be a multiset contained in $\Fpp \setminus \{\orig\}$ with $|A| = p+k$. If $\#\Sigma A < (k+2)p$, then there exists a nontrivial subgroup of $\Fpp$ that contains at least $p$ points of $A$, counting multiplicity.
\end{corollary}

Our methods of proof stem from two main ideas. First, we will obviously exploit the structure of $\Fpp$ as a direct sum of cyclic groups of prime order, within which we can apply the known Lemma~\ref{CD lemma} after using projections. Section~\ref{direct projects section} contains several elementary lemmas in this vein (see in particular Lemma~\ref{sweep to the left lemma}). It is important for us to utilize the flexibility coming from the fact that $\Fpp$ can be decomposed as the direct sum of two subgroups in many different ways. Second, our methods work best when there exists a single subgroup that contains many elements of the given multiset; however, by selectively replacing pairs of elements with their sums, we can increase the number of elements in a subgroup in a way that improves our lower bounds upon the sumset (see Lemma~\ref {any j lemma}). These methods, which appear in Section~\ref {thm proofs section}, combine to provide the proofs of Theorems \ref {conjecture true for k tiny} and~\ref{conjecture true for k small in terms of p}. Finally, Section~\ref{conjecture section} contains a generalization of Conjecture~\ref{2d conjecture} to higher-dimensional direct sums of $\F_p$, together with examples demonstrating that the conjecture would be best possible.

\section{Sumsets in abelian groups and direct products}
\label{direct projects section}

All of the results in this section are valid for general finite abelian groups and have correspondingly elementary proofs, although the last two lemmas seem rather less standard than the first few. In this section, $\GG$, $\HH$, and $\KK$ denote finite abelian groups, and $e$ denotes a group's identity element.

\begin{lemma}
\label{a la carte lemma}
Let $B_0,B_1,B_2,\dots,B_j$ be multisets in $\GG$, and set $A = B_0 \cup B_1 \cup \dots \cup B_j$. For each $1\le i\le j$, specify an element $x_i \in \Sigma B_i$, and set $C = B_0 \cup \{ x_1, \dots ,x_j\}$. Then $\Sigma C \subseteq \Sigma A$.
\end{lemma}

\begin{proof}
For each $1\le i\le j$, choose a submultiset $D_i \subseteq B_i$ such that the sum of the elements of $D_i$ equals $x_i$. By definition, every element $y$ of $\Sigma C$ equals the sum of the elements of some subset $E$ of $B_0$, plus $\sum_{i\in I} x_i$ for some $I\subseteq \{1,\dots,j\}$. But then $y$ equals the sum of the elements of $E \cup \bigcup_{i\in I} D_i$, which is an element of $\Sigma A$ since $E \cup \bigcup_{i\in I} D_i \subseteq B_0 \cup \bigcup_{1\le i\le j} B_i = A$.
\end{proof}

\begin{lemma}
\label{kneser bound}
Let $A_1,A_2,\dots,A_j$ be multisets in $\GG$, and set $A = A_1 \cup \dots \cup A_j$.  If $m$ is the largest cardinality of a proper subgroup of $\GG$, then either $\Sigma A = \GG$ or $\#\Sigma A \ge (\sum_{i=1}^j \# \Sigma A_i) - (j-1)m$.
\end{lemma}

\begin{proof}
Since $\Sigma A = \Sigma A_1 + \Sigma A_2 + \cdots + \Sigma A_j$ (viewed as ordinary sets), this follows immediately by inductive application of Kneser's theorem.
\end{proof}

For the remainder of this section, we will be dealing with groups that can be decomposed into a direct sum.

\begin{definition}
A subgroup $\HH$ of $\GG$ is called an {\em internal direct summand} if there exists a subgroup $\KK$ of $\GG$ such that $\GG$ is the internal direct sum of $\HH$ and $\KK$, or in other words, such that $\HH \cap \KK = \{e\}$ and $\HH + \KK = \GG$. Equivalently, $\HH$ is an internal direct summand of $\GG$ if there exists a {\em projection homomorphism} $\pi_\HH\colon \GG \to \HH$ that is the identity on $\HH$. Note that this projection homorphism does depend on the choice of $\KK$ but is uniquely determined by $\pi_\HH^{-1}(e) = \KK$.
\end{definition}

\begin{lemma}
\label{pi commutes with Sigma lemma}
For any homomorphism $f\colon \GG\to \HH$, and any subset $X$ of $\GG$, we have $f(\Sigma X) = \Sigma (f(X))$. In particular, if $\HH$ is an internal direct summand of $\GG$, then $\pi_\HH(\Sigma X) = \Sigma(\pi_\HH(X))$ for any subset $X$ of $\GG$.
\end{lemma}

\begin{proof}
Given $y\in f(\Sigma X)$, there exists $x\in \Sigma X$ such that $f(x)=y$. Hence we can find $x_1,\dots,x_j\in X$ such that $x_1+\cdots+x_j = x$, and so $f(x_1+\cdots+x_j) = y$. But $f$ is a homomorphism, and so $f(x_1)+\cdots+f(x_j) = y$, so that $y \in \Sigma(f(X))$. This shows that $f(\Sigma X) \subseteq \Sigma(f(X))$; the proof of the reverse inclusion is similar.
\end{proof}

\begin{lemma}
\label{pi split lemma}
Let $\GG = \HH \oplus \KK$, and let $D$ and $E$ be multisets contained in $\HH$ and $\KK$, respectively. For any $z\in \GG$,
\[
z \in \Sigma(D\cup E) \quad\text{if and only if}\quad \pi_\HH(z) \in\Sigma D \text{ and } \pi_\KK(z) \in\Sigma E.
\]
\end{lemma}

\begin{proof}
Since $z = \pi_\HH(z) + \pi_\KK(z)$, the ``if'' direction is obvious. For the converse, note that
\[
\pi_\HH(z) \in \pi_\HH\big( \Sigma (D\cup E) \big)  = \Sigma\big( \pi_\HH(D \cup E) \big)
\]
by Lemma~\ref {pi commutes with Sigma lemma}. On the other hand, $\pi_\HH(D) = D$ and $\pi_\HH(E) = \{e\}$, and so
\[
\pi_\HH(z) \in \Sigma\big( \pi_\HH(D) \cup \pi_\HH(E) \big) = \Sigma \big( D \cup \{e\} \big) = \Sigma D
\]
(since the sumset is not affected by whether $e$ is an allowed summand). A similar argument shows that $\pi_\KK(z) \in \Sigma E$, which completes the proof of the lemma.
\end{proof}

\begin{lemma}
\label{direct product lemma}
Let $\HH$ and $\KK$ be subgroups of $\GG$ satisfying $\HH \cap \KK = \{e\}$. Let $D$ and $E$ be multisets contained in $\HH$ and $\KK$, respectively. Then $\#\Sigma(D\cup E) = \#\Sigma D \cdot \#\Sigma E$.
\end{lemma}

\begin{proof}
Notice that every element of $\Sigma(D\cup E)$ is contained in $\HH+\KK$; therefore we may assume without loss of generality that $\GG = \HH \oplus \KK$. In particular, we may assume that $\HH$ and $\KK$ are internal direct summands of $\GG$, so that the projection maps $\pi_\HH$ and $\pi_\KK$ exist and every element $z\in \GG$ has a unique representation $z=x+y$ where $x\in \HH$ and $y\in \KK$; note that $x=\pi_\HH(z)$ and $y=\pi_\KK(z)$ in this representation.

To establish the lemma, it therefore suffices to show that $z = \pi_\HH(z) + \pi_\KK(z) \in \Sigma(D\cup E)$ if and only if $\pi_\HH(z) \in\Sigma D \text{ and } \pi_\KK(z) \in\Sigma E$; but this is exactly the statement of Lemma~\ref{pi split lemma}.
\end{proof}

The next lemma is a bit less standard yet still straightforward: in a direct product of two abelian groups, it characterizes the elements of a sumset that lie in a given coset of one of the direct summands.

\begin{lemma}
\label{orthogonal structure lemma}
Let $\HH$ and $\KK$ be subgroups of $\GG$ satisfying $\HH \cap \KK = \{e\}$. Let $D$ and $E$ be multisets contained in $\HH$ and $\KK$, respectively. For any $y\in \KK$:
\begin{enumerate}
 \item if $y\in \Sigma E$, then $(\HH+\{y\}) \cap \Sigma(D\cup E) = \Sigma D + \{y\}$;
 \item if $y\notin \Sigma E$, then $(\HH+\{y\}) \cap \Sigma(D\cup E) = \emptyset$.
\end{enumerate}
\end{lemma}

\begin{proof}
As in the proof of Lemma~\ref {direct product lemma}, we may assume without loss of generality that $\GG = \HH \oplus \KK$. Suppose that $z$ is an element of $(\HH+\{y\}) \cap \Sigma(D\cup E)$. Since $z\in \HH+\{y\}$, we may write $z=x+y$ for some $x\in\HH$, whence $\pi_\KK(z) = \pi_\KK(x) + \pi_\KK(y) = e+y = y$. On the other hand, since $z\in \Sigma(D\cup E)$, we see that $y \in \Sigma E$ by Lemma~\ref{pi split lemma}. In other words, the presence of any element $z\in (\HH+\{y\}) \cap \Sigma(D\cup E)$ forces $y\in \Sigma E$, which establishes part (b) of the lemma.

We continue under the assumption $y\in \Sigma E$ to prove part (a). The inclusions $\Sigma D + \{y\} \subseteq \HH+\{y\}$ and $\Sigma D + \{y\} \subseteq \Sigma(D\cup E)$ are both obvious, and so $\Sigma D + \{y\} \subseteq (\HH+\{y\}) \cap \Sigma(D\cup E)$. As~for the reverse inclusion, let $z \in (\HH+\{y\}) \cap \Sigma(D\cup E)$ as above; then $\pi_\HH(z) \in \Sigma D$ by Lemma~\ref{pi split lemma}, whence $z = \pi_\HH(z) + \pi_\KK(z) = \pi_\HH(z) + y \in \Sigma D + \{y\}$ as required.
\end{proof}

Finally we can establish the lemma that we will make the most use of when we return to the setting $\GG=\Fpp$ in the next section.

\begin{lemma}
\label{sweep to the left lemma}
Let $\GG = \HH \oplus \KK$, and let $C$ be a multiset contained in $\GG$. Let $D=C\cap \HH$, let $F = C \setminus D$, and let $E = \pi_\KK(F)$. Then $\#\Sigma C \ge \#\Sigma D \cdot \#\Sigma E$.
\end{lemma}

\begin{proof}
Lemma~\ref {direct product lemma} tells us that $\#\Sigma (D\cup E) = \#\Sigma D \cdot \#\Sigma E$, and so it suffices to show that $\#\Sigma C \ge \#\Sigma (D\cup E)$. We accomplish this by showing that
\begin{equation}
\# \big( (\HH + \{y\}) \cap \Sigma C \big) \ge \# \big( (\HH + \{y\}) \cap \Sigma (D\cup E) \big)
\label{one line at a time}
\end{equation}
for all $y\in\KK$.

For any $y\in\KK \setminus \Sigma E$, Lemma~\ref{orthogonal structure lemma} tells us that $(\HH + \{y\}) \cap \Sigma (D\cup E) = \emptyset$, in which case the inequality~\eqref{one line at a time} holds trivially. For any $y \in \Sigma E$, Lemma~\ref{orthogonal structure lemma} tells us that $(\HH + \{y\}) \cap \Sigma (D\cup E) = \Sigma D + \{y\}$, and so the right-hand side of the inequality~\eqref{one line at a time} equals $\#\Sigma D$.

On the other hand, since $\Sigma E = \Sigma (\pi_\KK(F)) = \pi_\KK(\Sigma F)$ by Lemma~\ref{pi commutes with Sigma lemma}, there exists at least one element $z\in \Sigma F$ satisfying $\pi_\KK(z)=y$; as $\GG = \HH \oplus \KK$, this is equivalent to saying that $z \in \HH + \{y\}$. Since $\Sigma D \subseteq \HH$, we have $\Sigma D + \{z\} \subseteq \HH + \{y\}$ as well. But the inclusion $\Sigma D + \{z\} \subseteq \Sigma D + \Sigma F = \Sigma C$ is trivial, and therefore $\Sigma D + \{z\} \subseteq (\HH + \{y\}) \cap \Sigma C$; in particular, the left-hand side of the inequality~\eqref{one line at a time} is at least $\#\Sigma D$. Combined with the observation that the right-hand side equals $\#\Sigma D$, this lower bound establishes the inequality~\eqref{one line at a time} and hence the lemma.
\end{proof}

These lemmas might be valuable for studying sumsets in more general abelian groups. They will prove to be particularly useful for studying sumsets in $\Fpp$, however, essentially because there are many ways of writing $\Fpp$ as an internal direct sum of two subgroups (which are simply lines through $\orig$).

\section{Lower bounds for sumsets}
\label{thm proofs section}

In this section we establish Theorems~\ref {conjecture true for k tiny} and~\ref {conjecture true for k small in terms of p}; the proofs employ two combinatorial propositions which we defer to the next section. It would be possible to prove these two theorems at the same time, at the expense of a bit of clarity; however, we find it illuminating to give complete proofs of Theorem~\ref {conjecture true for k tiny} (the cases $|A|=p$ and $|A|=p+1$) first, as the proofs will illustrate the methods used to prove the more general Theorem~\ref {conjecture true for k small in terms of p}. Seeing the limitations of the proof of Theorem~\ref {conjecture true for k tiny} will also motivate the formulation of our main technical tool, Lemma~\ref {any j lemma}.

Throughout this section, $A$ will denote a valid multiset contained in $\Fpp$. For any $x\in\Fpp$, we let $\lin x$ denotes the subgroup of $\Fpp$ generated by $x$ (that is, the line passing through both the origin $\orig$ and $x$), and we let $m_x$ denote the multiplicity with which $x$ appears in $A$, so that $|A| = \sum_{x\in\Fpp} m_x$. The fact that $A$ is valid means that $m_\orig=0$ and $\sum_{t\in\lin x} m_t < p$ for every $x\in\Fpp\noto$.

Our first lemma quantifies the notion that we can establish sufficiently good lower bounds for the cardinality of $\Sigma A$ if we know that there are enough elements of $A$ lying in one subgroup of~$\Fpp$. Naturally, the method of proof is to partition $A$ into the elements lying in that subgroup and all remaining elements, project the remaining elements onto a complementary subgroup, and then use Lemma~\ref{CD lemma} in each subgroup separately.

\begin{lemma}
\label{conjecture true if enough on line lemma}
Let $A$ be any valid multiset contained in $\Fpp$.
Suppose that for some $x\in\Fpp\noto$,
\begin{equation}
\label{symmetric bound}
\sum_{y\in\lin x} m_y \ge |A| - (p-1).
\end{equation}
Then $\#\Sigma A \ge (|A|+2-p)p$.
\end{lemma}

\begin{remark}
The conclusion is worse than trivial if $|A| < p-1$; also, the fact that $A$ is valid means that the left-hand side of equation~\eqref{symmetric bound} is at most $p-1$, and so the lemma is vacuous if $|A| > 2p-2$. Therefore in practice the lemma will be applied only to multisets $A$ satisfying $p-1 \le |A| \le 2p-2$.
\end{remark}

\begin{proof}
Let $D = A \cap \lin x$; note that $|D| \le p-1$ since $A$ is a valid multiset, and note also that $|D| = \sum_{y\in\lin x} m_y \ge |A| - (p-1)$ by assumption. Set $F = A \setminus D$. Choose any nontrivial subgroup $\KK$ of $\Fpp$ other than $\lin x$, and set $E = \pi_\KK(F)$. Then by Lemma~\ref{sweep to the left lemma}, we know that $\#\Sigma A \ge \#\Sigma D \cdot \#\Sigma E$. By Lemma~\ref {CD lemma} and the fact that $\orig\notin D\cup E$, we obtain
\begin{align}
\#\Sigma A &\ge \min \big\{ p, 1 + |D| \big\} \cdot \min \big\{ p, 1 + |E| \big\} \notag \\
&= \min \big\{ p, 1 + |D| \big\} \cdot \min \big\{ p, 1 + |A| - |D| \big\},
\label{actually special case}
\end{align}
since $|E| = |F| = |A|-|D|$. The inequalities $|D| \le p-1$ and $|A| - |D| \le p-1$ ensure that $p$ is the larger element in both minima, and so we have simply
\[
\#\Sigma A \ge (1+|D|)(1+|A|-|D|) = \tfrac14|A|^2 + |A| + 1 - \big( |D| - \tfrac12|A| \big)^2.
\]
The pair of inequalities $|D| \le p-1$ and $|A| - |D| \le p-1$ is equivalent to the inequality $\big| |D| - \tfrac12|A| \big| \le p-1- \tfrac12|A|$; therefore
\[
|\Sigma A| \ge \tfrac14|A|^2 + |A| + 1 - \big( p-1- \tfrac12|A| \big)^2 = (|A|+2-p)p,
\]
as claimed.
\end{proof}

This lemma alone is sufficient to establish Theorem~\ref{conjecture true for k tiny}.

\begin{proof}[Proof of Theorem~\ref{conjecture true for k tiny}(a)]
When $|A|=p$, the right-hand side of the inequality~\eqref {symmetric bound} equals 1, and so the inequality holds for any $x\in A$. Therefore Lemma~\ref {conjecture true if enough on line lemma} automatically applies, yielding $\#\Sigma A \ge (|A|+2-p)p = 2p$ as desired. (In fact essentially the same proof gives the more general statement: if $A$ is a multiset contained in $\Fpp\setminus\{\orig\}$ but not contained in any proper subgroup, and $|A|\ge p$, then $\#\Sigma A \ge 2|A|$.)
\end{proof}

\begin{proof}[Proof of Theorem~\ref{conjecture true for k tiny}(b)]
We are assuming that $|A| = p+1$. Suppose first that there exists a nontrivial subgroup of $\Fpp$ that contains at least two points of $A$ (including possibly two copies of the same point). Choosing any nonzero element $x$ in that subgroup, we see that the inequality~\eqref {symmetric bound} is satisfied, and so Lemma~\ref {conjecture true if enough on line lemma} yields $\#\Sigma A \ge (|A|+2-p)p = 3p$ as desired.

From now on we may assume that there does not exist a nontrivial subgroup of $\Fpp$ that contains at least two points of $A$. Since there are only $p+1$ nontrivial subgroups of $\Fpp$, it must be the case that $A$ consists of exactly one point from each of these $p+1$ subgroups; in particular, the elements of $A$ are distinct. We can verify the assertion for $p\le11$ by exhaustive computation (see the remarks after the end of this proof), so from now on we may assume that $p\ge13$.

Suppose first that all sums of pairs of distinct elements from $A$ are distinct. All these sums are elements of $\Sigma A$, and thus $\#\Sigma A \ge \binom{p+1}2 >3p$ since $p\ge13$.

The only remaining case is when two pairs of distinct elements from $A$ sum to the same point of $\Fpp$. Specifically, suppose that there exist $x_1,y_1,x_2,y_2\in A$ such that $x_1+y_1=x_2+y_2$. Partition $A = B_0 \cup B_1 \cup B_2$ where $B_1=\{x_1,y_1\}$ and $B_2=\{x_2,y_2\}$ and hence $B_0 = A \setminus \{x_1,y_1,x_2,y_2\}$; note that this really is a partition of $A$, as the fact that $x_1+y_1=x_2+y_2$ forces all four elements to be distinct. Moreover, if we define $z=x_1+y_1=x_2+y_2$, then we know that $z\ne\orig$ since $x_1$ and $y_1$ are in different subgroups.

Define $C$ to be the multiset $B_0 \cup \{ z, z \}$; by Lemma~\ref {a la carte lemma}, we know that $\#\Sigma A \ge \#\Sigma C$. Define $D = C \cap \lin z$; we claim that $|D| = 3$. To see this, note that $A$ has exactly one point in every nontrivial subgroup, and in particular $A$ has exactly one point in $\lin z$. Furthermore, that point cannot be $x_1$ for example, since then $y_1 = z-x_1$ would also be in that subgroup; similarly that point cannot be $x_2$, $y_1$, or $y_2$. We conclude that $B_0$ has exactly one point in $\lin z$, whence $C$ has exactly three points in $\lin z$.

Now define $F = C \setminus D$, so that $|F| = |C| - |D| = (|B_0|+2)-3 = (|A|-4+2)-3 = p-4$. Let $\KK$ be any nontrivial subgroup other than $\lin z$, and set $E = \pi_\KK(F)$. The lower bounds $\#\Sigma D \ge 4$ and $\#\Sigma E \ge p-3$ then follow from Lemma~\ref {CD lemma}. By Lemma~\ref{sweep to the left lemma}, we conclude that $\#\Sigma C \ge \#\Sigma D \cdot \#\Sigma E = 4(p-3) > 3p$ since $p\ge13$.
\end{proof}

\begin{remark}
The computation that verifies Theorem~\ref{conjecture true for k tiny}(b) for $p\le11$ should be done a little bit intelligently, since there are $10^{12}$ subsets $A$ of $\F_{11}^2$ (for example) consisting of exactly one nonzero element from each nontrivial subgroup. We describe the computation in the hardest case $p=11$. Let us write the elements of $\F_{11}^2$ as ordered pairs $(s,t)$ with $s$ and $t$ considered modulo~$11$. By separately dilating the two coordinates of $\F_{11}^2$ (which does not alter the cardinality of $\Sigma A$), we may assume without loss of generality that $A$ contains both $(1,0)$ and $(0,1)$. We also know every such $A$ contains a subset of the form $\{ (i,i), (j,2j), (k,3k), (\ell,4\ell) \}$ for some integers $1\le i,j,k,\ell \le 10$. Therefore the cardinality of every such $\Sigma A$ is at least as large as the cardinality of one of the subsumsets $\Sigma \big( \{ (1,0), (0,1), (i,i), (j,2j), (k,3k), (\ell,4\ell) \} \big)$.

There are $10^4$ such subsumsets, and direct computation shows that all of them have more than $33$ distinct elements except for the sixteen cases $\Sigma \big( \{ (1,0), (0,1), \pm(1,1), \pm(1,2), \pm(1,3), \pm(1,4) \} \big)$, which each contain $32$ distinct elements. It is then easily checked that any subsumset of the form $\Sigma \big( \{ (1,0), (0,1), \pm(1,1), \pm(1,2), \pm(1,3), \pm(1,4), (m,5m) \} \big)$ with $1\le m\le10$ contains more than 33 distinct elements. This concludes the verification of Theorem~\ref{conjecture true for k tiny}(b) for $p=11$, and the cases $p\le7$ are verified even more quickly.
\end{remark}

We now foreshadow the proof of Theorem~\ref{conjecture true for k small in terms of p} by reviewing the structure of the proof of Theorem~\ref{conjecture true for k tiny}(b). In that proof, we quickly showed that the desired lower bound held if there were enough elements of $A$ in the same subgroup. Also, the desired lower bound certainly held if there were enough distinct sums of pairs of elements of~$A$. If however no subgroup contained enough elements of $A$ and there were only a few distinct sums of pairs of elements of $A$, then we showed that we could find multiple pairs of elements summing to the same point in $\Fpp$. Replacing those elements in $A$ with multiple copies of their joint sum, we found that the corresponding subgroup now contained enough elements to carry the argument through.

The following lemma quantifies the final part of this strategy, where we replace $j$ pairs of elements of $A$ with their joint sum and then use our earlier ideas to bound the cardinality of the sumset from below.

\begin{lemma}
\label{any j lemma}
Let $A$ be any valid multiset contained in $\Fpp$, and let $z\in\Fpp\noto$. For any integer $j$ satisfying
\begin{equation}
0\le j \le \tfrac12 \sum_{t\in\Fpp\setminus\lin z} \min\{m_t,m_{z-t}\},
\label{allows choosing pairs}
\end{equation}
we have
\[
\#\Sigma A \ge \min\bigg\{ p, 1 + j + \sum_{y\in\lin z} m_y \bigg\} \min\bigg\{ p, 1 + |A| - 2j - \sum_{y\in\lin z} m_y \bigg\}.
\]
\end{lemma}

\begin{remark}
This can be seen as a generalization of Lemma~\ref {conjecture true if enough on line lemma}, as equation~\eqref {actually special case} is the special case $j=0$ of this lemma.
\end{remark}

\begin{proof}
Partition $A = B_0 \cup B_1 \cup \cdots \cup B_j$, where for each $1\le i\le j$, the multiset $B_i$ has exactly two elements, neither contained in $\lin z$, that sum to $z$ (the complimentary submultiset $B_0$ is unrestricted). The upper bound~\eqref{allows choosing pairs} for $j$ is exactly what is required for such a partition to be possible; the factor of $\frac12$ arises because the sum on the right-hand side of~\eqref{allows choosing pairs} double-counts the pairs $(t,z-t)$ and $(z-t,t)$. Then set $C$ equal to $B_0$ with $j$ additional copies of $z$ inserted. By Lemma~\ref{a la carte lemma}, we know that $\#\Sigma A \ge \#\Sigma C$.

Now let $D$ be the intersection of $C$ with the subgroup $\lin z$, and let $F = C \setminus D$. Let $\KK$ be any nontrivial subgroup other than $\lin z$, and set $E = \pi_\KK(F)$. By Lemma~\ref{sweep to the left lemma}, we know that $\#\Sigma C \ge \#\Sigma D \cdot \#\Sigma E$. However, the number of elements of $D$ (counting multiplicity) is $j$ more than the number of elements of $B_0 \cap \lin z$; this is the same as $j$ more than the number of elements of $A \cap \lin z$ (since no elements of $B_1,\dots,B_j$ lie on $\lin z$), or in other words $j + \sum_{y\in\lin z} m_y$. Similarly, the number of elements of $E$ (equivalently, of $F$) is equal to the number of elements of $B_0 \setminus \lin z$; this is the same as $2j$ less than the number of elements of $A \setminus \lin z$, or in other words $|A| - 2j - \sum_{y\in\lin z} m_y$. The lower bounds $\#\Sigma D \ge \min\big\{ p, 1 + j + \sum_{y\in\lin z} m_y \big\}$ and $\#\Sigma E \ge \min\big\{ p, 1 + |A| - 2j - \sum_{y\in\lin z} m_y \big\}$ then follow from Lemma~\ref {CD lemma}; the chain of inequalities $\#\Sigma A \ge \#\Sigma C \ge \#\Sigma D \cdot \#\Sigma E$ establishes the lemma.
\end{proof}

We are now ready to use Lemma~\ref{any j lemma} to establish Conjecture~\ref{2d conjecture} when $|A|=p+k$, for all but finitely many primes $p$ depending on~$k$. Let $H_k = 1 + \tfrac12 + \cdots + \tfrac1k$ denote the $k$th harmonic number.

\begin{theorem}
\label{conjecture true for p large in terms of k}
Let $k\ge2$ be any integer, and let $A$ be any valid multiset contained in $\Fpp$ such that $|A| = p+k$. If $p\ge 4(k+1)^2 H_k - 2k$, then $\#\Sigma A \ge (k+2)p$.
\end{theorem}

\begin{remark}

Using the elementary bound $H_k \le \gamma + \log(k+1)$, where $\gamma$ denotes the Euler--Mascheroni constant, we see that Theorem~\ref{conjecture true for p large in terms of k} holds as long as $p \ge 4(k+1)^2 (\gamma + \log(k+1))$. Theorem~\ref{conjecture true for k small in terms of p} can thus be readily deduced from Theorem~\ref {conjecture true for p large in terms of k} as follows: If $k+1 \le \sqrt{p/(2\log p+1)}$ then $p \ge 4(k+1)^2 (\tfrac14 + \tfrac12 \log p)$.  In this case we certainly have $p \ge (1 + 2\log 2) (k+1)^2$, whence $\log p \ge \frac45 + 2 \log(k+1)$ and $\tfrac14 + \tfrac12 \log p \ge \gamma + \log(k+1)$.

\end{remark}

\begin{proof}
If there are $k+1$ elements of $A$ in some nontrivial subgroup, then we are done by Lemma~\ref {conjecture true if enough on line lemma}. Therefore we may assume that there are at most $k$ points in each subgroup; in particular, $m_x\le k$ for all $x\in\Fpp$. We now argue that if $\Sigma A$ is small, then there must be lots of pairs of elements of $A$ that add to the same element of $\Fpp$, at which point we will be able to invoke Lemma~\ref {any j lemma}. We may assume that $\Sigma A \ne \Fpp$, for otherwise we are done immediately.

For each $1 \le i \le k$, we define the level set $A_i = \{x \in \Fpp : m_x \ge i\}$. Notice that $A$ can be written precisely as the multiset union $A_1 \cup A_2 \cup \cdots \cup A_k$, and so $\sum_{i=1}^k \#A_i = |A| = p+k$. Let $B_i$ be the multiset formed by the sums of pairs of elements of $A_i$ not in the same subgroup:
\[
B_i = \big\{ x + y \colon x,y\in A_i,\, \lin x\ne\lin y \big\}.
\]
Note that $\orig \notin B_i$ (the restriction $\lin x\ne\lin y$ ensures that $x \ne -y$) and that every element of $B_i$ occurs with even multiplicity (the restriction $\lin x\ne\lin y$ ensures that $x \ne y$).  It is not hard to estimate the relative sizes of $\#A_i$ and $|B_i|$: for each $x \in A_i$ there are at most $k$ elements of $A$ lying in the subgroup $\lin x$.  Since each such $x$ occurs with multiplicity at least $i$ in $A$, there are at most $k/i$ distinct values of $y$ excluded by the condition $\lin x \ne \lin y$.  Hence
$|B_i| \ge \#A_i (\#A_i - k/i)$, which implies that
\begin{equation}\label{bound on a_i}
\#A_i \le \frac{k}{i} + \sqrt{|B_i|}.
\end{equation}
Since $\sum_{i=1}^k \#A_i$ is fixed, this shows that $|B_i|$ cannot be very small on average.  At the same time, $\#B_i$ cannot get very large: if $\sum_{i=1}^k \#B_i \ge (2k+1)p$, then (under our assumption that $\Sigma A \ne \Fpp$) Lemma~\ref{kneser bound} already yields
\[
\#\Sigma A \ge \sum_{i=1}^k \#\Sigma A_i - (k-1)p > \sum_{i=1}^k \#B_i - (k-1)p \ge (k+2)p.
\]
where the middle inequality holds because $B_i \subseteq \Sigma A_i$. We may therefore assume henceforth that
\begin{equation} \label{bound on b_i}
\sum_{i=1}^k \#B_i < (2k+1)p.
\end{equation}

Let us now introduce the weighted height parameter
\begin{equation} \label{eta def}
\eta = \max_{1\le i \le k} \left\{ \frac{i|B_i|}{2\#B_i} : \#B_i > 0 \right\}.
\end{equation}
We shall show shortly that $\eta > k+1$.  Assuming so, then for some $1 \le i \le k$, we have
\[
\frac{|B_i|}{2\#B_i} > \frac{k+1}{i},
\]
so by the pigeonhole principle, there exists some $z \in B_i$ (in particular $z\ne\orig$) occurring with multiplicity greater than $2(k+1)/i$; since this multiplicity is an even integer, it must be at least $2(k+2)/i.$  For each solution $x+y = z$ corresponding to an occurrence of $z$ in $B_i$, we have by construction that $x,y \notin \lin z$ and $m_x, m_y \ge i$, so for this particular choice of~$z$,
\[
\tfrac12 \sum_{t \in \Fpp \setminus\lin z } \min\{m_t, m_{z-t}\} \ge k+2.
\]
Furthermore, $\sum_{y\in\lin z} m_y \le k$ by assumption. Therefore we are free to apply Lemma~\ref{any j lemma} with $j = (k+2) - \sum_{y\in\lin z} m_y,$
which gives the lower bound
\[
\#\Sigma A \ge \min\{p,k+3\} \min\bigg\{p,p - k - 3  + \sum_{y\in\lin z} m_y\bigg\} \ge (k+3)(p-k-3) \ge (k+2)p
\]
(the last step used the inequality $p \ge (k+3)^2$, which certainly holds under the hypotheses of the theorem).

It remains only to verify that $\eta > k+1$. Summing the inequality~\eqref{bound on a_i} over all $1\le i\le k$ yields
\[
p+k = \sum_{i=1}^k \#A_i \le k H_k + \sum_{i=1}^k \sqrt{|B_i|} \le kH_k + \sqrt{2\eta} \sum_{i=1}^k \sqrt{\frac{\#B_i}{i}},
\]
using the definition~\eqref{eta def} of~$\eta$. We estimate the rightmost sum using Cauchy--Schwarz together with the inequality~\eqref{bound on b_i}:
\[
\sum_{i=1}^k \sqrt{\frac{\#B_i}{i}} \le \bigg(\sum_{i=1}^k \#B_i\bigg)^{1/2}\bigg(\sum_{i=1}^k \frac1i\bigg)^{1/2} < \sqrt{(2k+1) p H_k}.
\]
Combining the previous two inequalities gives $p+k - kH_k < \sqrt{ \eta(4k+2) p H_k}$, so that
\[
\eta > \frac{(p+k-kH_k)^2}{(4k+2) pH_k} > \frac{p(p+2(k-kH_k))}{(4k+2) pH_k} = \frac{(p+2k)-2kH_k}{(4k+2)H_k} \ge \frac{4(k+1)^2H_k-2kH_k}{(4k+2)H_k}
\]
by the hypothesis on the size of~$p$. In other words,
\[
\eta > \frac{2(k+1)^2-k}{2k+1} = k+1+\frac1{2k+1},
\]
which completes the proof of the theorem.
\end{proof}

\section{A wider conjecture}
\label{conjecture section}

As mentioned earlier, Conjecture~\ref{2d conjecture} is just one part of a more far-reaching conjecture concerning sumsets of multisets in $\F_p^m$. Before formulating that wider conjecture, we must expand the definition of a valid multiset to $\F_p^m$.

\begin{definition}
\label{more valid definition}
Let $p$ be an odd prime, and let $m$ be a positive integer. A multiset $A$ contained in $\F_p^m$ is {\em valid} if:
\begin{itemize}
\item $\orig\notin A$; and
\item for each $1\le d\le m$, every subgroup of $\F_p^m$ that is isomorphic to $\F_p^d$ contains fewer than $dp$ points of $A$, counting multiplicity.
\end{itemize}
\end{definition}

\noindent When $m=1$, a multiset contained in $\F_p$ is valid precisely when it does not contain $0$; when $m=2$ and $|A| < 2p$, this definition of valid agrees with Definition~\ref{valid definition} for multisets contained in $\Fpp$. Note that in particular, Definition~\ref{more valid definition}(b) implies that every valid multiset contained in $\F_p^m$ has at most $mp-1$ elements, counting multiplicity. We now give an example showing that this upper bound $mp-1$ can in fact be achieved. Throughout this section, let $\{x_1,\dots,x_m\}$ denote a generating set for $\F_p^m$, and let $\KK_d = \lin{x_1,\dots,x_d}$ denote the subgroup of $\F_p^m$ generated by $\{x_1,\dots,x_d\}$, so that $\KK_d \cong \F_p^d$.

\begin{example}
\label{valid example}
Let $A_1$ be the multiset consisting of $p-1$ copies of $x_1$; for $2\le j\le m$ let $A_j = \{ x_j + ax_1 \colon 0\le a\le p-1 \}$; and define $B_m = \bigcup_{j=1}^m A_j$. Then $|B_m| = (p-1) + (m-1)p = mp-1$ and $\orig\notin B_m$. To verify that $B_m$ is a valid subset of $\F_p^m$, let $\HH$ be any subgroup of $\F_p^m$ that is isomorphic to $\F_p^d$; we need to show that $B_m$ contains fewer than $dp$ points of $\HH$.

First suppose that $x_1\notin\HH$, which implies that $bx_1\notin\HH$ for every nonzero multiple $bx_1$ of~$x_1$. Then for each $2\le j\le m$, at most one of the elements of $A_j$ can be in $\HH$, since the difference of any two such elements is a nonzero multiple of $x_1$. Therefore $|B_m\cap\HH| = \ell$ for some $1\le\ell\le m-1$, and in fact all $\ell$ of these elements are of the form $x_j+ax_1$ for $\ell$ distinct values of~$j$. Since no such element is in the subgroup spanned by the others, we conclude that $d\ge\ell$, and so the necessary inequality $|B_m\cap\HH|=\ell\le d<dp$ is amply satisfied.

Now suppose that $x_1\in\HH$. Then for each $2\le j\le m$, either all or none of the elements of $A_j$ are in $\HH$. By reindexing the $x_i$, we may choose an integer $1\le\ell\le m$ such that $\HH$ contains $A_1 \cup \cdots \cup A_\ell$ and is disjoint from $A_{\ell+1} \cup \cdots \cup A_m$. In particular, $|B_m\cap\HH| = (p-1)+(\ell-1)p = \ell p-1$. But $\HH$ contains $\{x_1,\dots,x_\ell\}$ and hence $d\ge\ell$, so that $\ell p-1\le dp-1$ as required.
\end{example}

We may now state our wider conjecture; Conjecture~\ref{2d conjecture} is the special case $q=1$ of part (a) of this conjecture.

\begin{conjecture}
\label {any d conjecture}
Let $p$ be an odd prime. Let $m$ be a positive integer, and let $A$ be a valid multiset of $\F_p^m$ with $|A|\ge p$. Write $|A| = qp+k$ with $0\le k\le p-1$.
\begin{enumerate}
\item If $0\le k\le p-3$, then $\#\Sigma A \ge (k+2)p^q$.
\item If $k=p-2$, then $\#\Sigma A \ge p^{q+1}-1$.
\item If $k=p-1$, then $\#\Sigma A \ge p^{q+1}$.
\end{enumerate}
In particular, if $|A| = mp-1$ then $\Sigma A = \F_p^m$.
\end{conjecture}

We remark that the quantity $dp$ in Definition~\ref{more valid definition}, bounding the number of elements in a valid multiset that can lie in a rank-$d$ subgroup, has been carefully chosen in light of this conjecture: by Conjecture~\ref {any d conjecture}(c), any valid multiset $A$ with at least $dp-1$ elements counting multiplicity must satisfy $\#\Sigma A \ge p^d$. In particular, if $A$ is a valid multiset contained in a subgroup $\HH < \F_p^m$ that is isomorphic to $\F_p^d$, then $|A|\ge dp-1$ implies that $\Sigma A = \HH$. Therefore allowing $dp$ elements in such a subgroup would always be ``wasteful''. Of course, the validity of Definition~\ref{more valid definition} for rank-$d$ subgroups depends crucially upon the truth of Conjecture~\ref {any d conjecture}(c) for $q=d-1$.

The conjecture is restricted to multisets $A$ with $|A|\ge p$ because we already know the truth for smaller multisets, for which the definition of ``valid'' is simply the condition that $\orig \notin A$: when $|A|\le p-1$, the best possible lower bound is $\#\Sigma A \ge |A|+1$ as in Lemma~\ref {CD lemma}. We remark that Peng~\cite[Theorem 2]{P1} has proved Conjecture~\ref {any d conjecture}(c) in the case $m=2$ and $q=1$, under even a slightly weaker hypothesis; in other words, he has shown that if $A$ is a valid multiset contained in $\Fpp$ with $|A| = 2p-1$, then $\Sigma A = \Fpp$. (We remark that Mann and Wou~\cite{MW} have proved in the case that $A$ is actually a set---that is, a multiset with distinct elements---that $\#A = 2p-2$ suffices to force $\Sigma A = \Fpp$.) Peng considers the higher-rank groups $\F_p^m$ as well, but the multisets he allows (see \cite[Theorem 1]{P2}) form a much wider class than our valid multisets, and so his conclusions are much weaker than Conjecture~\ref{any d conjecture} for $q\ge2$. Finally, we mention that we have completely verified Conjecture~\ref{any d conjecture} by exhaustive computation for the groups $\F_p^2$ with $p\le 7$ and also for the group $\F_3^3$.

It is easy to see that all of the lower bounds in Conjecture~\ref{any d conjecture}(a), if true, would be best possible. Given $q\ge1$ and $0\le k\le p-3$, let $A'$ be any valid multiset contained in $\KK_q$ with $|A'| = qp-1$ (such as the one given in Example~\ref{valid example} with $m=q$), and let $A$ be the union of $A'$ with $k+1$ copies of $x_{q+1}$. Then $\Sigma A = \{ y + ax_{q+1}\colon y\in \Sigma A',\, 0\le a\le k+1\}$ and thus $\#\Sigma A' = (k+2)\#\Sigma A \le (k+2)p^q$ since $\Sigma A$ is contained in~$\KK_q$. Similarly, the fact that there exists a valid multiset contained in $\KK_{q+1}$ with $qp+(q-1)=(q+1)p-1$ elements (such as the one given in Example~\ref{valid example} with $m=q+1$) shows that the lower bound in Conjecture~\ref{any d conjecture}(c) would be best possible, since the sumset of this multiset would still be contained in $\KK_{q+1}$ and thus would have at most $p^{q+1}$ distinct elements.

The lower bound in Conjecture~\ref{any d conjecture}(b) might seem counterintuitive, especially in comparison with the pattern established in Conjecture~\ref{any d conjecture}(a). However, we can give an explicit example showing that the lower bound $p^{q+1}-1$ for $\#\Sigma A$ cannot be increased:

\begin{example}
\label{border example}
When $p$ is an odd prime, define $B'_m$ to be the set $B_m$ from Example~\ref{valid example} with one copy of $x_1$ removed, so that $B'_m$ contains $x_1$ with multiplicity only $p-2$. Since $B_m$ is a valid multiset contained in $\F_p^m$, so is $B'_m$. We have $|B'_m| = |B_m|-1 = (mp-1)-1 = (m-1)p + (m-2)$, and we claim that $-x_1\notin \Sigma B'_m$; this will imply that $\#\Sigma B'_m \le p^m-1$, and so the lower bound for $\#\Sigma A$ in Conjecture~\ref{any d conjecture}(b) cannot be increased. (In fact it is not hard to show that every other element of $\F_p^m$ is in $\Sigma B'_m$, and so $\#\Sigma B'_m$ is exactly equal to $p^m-1$.)

Suppose for the sake of contradition that $-x_1\in \Sigma B'_m$, and let $C$ be a submultiset of $B'_m$ such that $-x_1 = \sum_{y\in C}y$. For each $2\le j\le m$, define $\ell_j = |C\cap A_j|= \#\big( C \cap \{ x_j + ax_1\colon 0\le a\le p-1 \} \big)$. Then
\[
-x_1 = \sum_{y\in C} y = t x_1 + \ell_2 x_2 + \ell_3 x_3 + \cdots + \ell_m x_m
\]
for some integer~$t$. It follows from this equation that each $\ell_j$ must equal either $0$ or $p$. However, if $\ell_j = p$ then
\[
\sum_{y \in C\cap A_j} y = \sum_{0\le a\le p-1} (x_j + ax_1) = px_j + \frac{p(p-1)}2 x_1 = \orig
\]
(since $p$ is odd). So in either case, if $s = |C\cap A_1|$ is the multiplicity with which $x_1$ appears in $C$, then
\[
-x_1 = \sum_{y\in C} y = s x_1 + \sum_{j=2}^m \sum_{y\in C\cap A_j} y = s x_1 + \orig + \cdots + \orig.
\]
This is a contradiction, however, since $s$ must lie between $0$ and $p-2$. Therefore $-x_1$ is indeed not an element of $\Sigma B'_m$, as claimed.
\end{example}

The line of questioning in this section turns out to be uninteresting when $p=2$: when the multiset $A$ does not contain $\orig$, the condition that no rank-$1$ subgroup of $\F_2^m$ contain $2$ points of $A$ is simply equivalent to $A$ not containing any element with multiplicity greater than~$1$. It is easy to check that if $A$ consists of any $q$ points in $\F_2^m$ that do not lie in any subgroup isomorphic to $\F_2^{q-1}$, then $\Sigma A$ fills out the entire rank-$q$ subgroup generated by~$A$. In other words, the analogous definition of ``valid'' for multisets in $\F_2^m$ would simply be a set of $q$ points that generate a rank-$q$ subgroup of $\F_2^m$, and we would always have $\#\Sigma A = 2^{|A|} = 2^{\#A}$ for valid (multi)sets in~$\F_2^m$.

\section*{Acknowledgments}

The collaboration leading to this paper was made possible thanks to Jean--Jacques Risler, Richard Kenyon, and especially Ivar Ekeland; the authors also thank the University of British Columbia and the Institut d'\'Etudes Politiques de Paris for their undergraduate exchange program. 
The first author thanks Andrew Granville for conversations that explored this topic and eventually led to the formulation of the conjectures herein.

\end{document}